\documentclass[11pt, oneside]{amsart}
\usepackage{latexsym}
\usepackage{enumitem}
\usepackage{layout}
\usepackage{amsfonts}
\usepackage{amssymb}
\usepackage{epsfig}
\usepackage{color}
\usepackage{cite}

\usepackage[normalem]{ulem}
\usepackage[perpage,para,symbol*]{footmisc}
\usepackage{fancyhdr}
\usepackage{times}
\usepackage{lipsum} 
\newcommand{\R}{\mathbf{R}}
\newcommand{\E}{\mathcal{E}}
\newcommand{\C}{\mathbb{C}}

\numberwithin{equation}{section}
\theoremstyle{plain}
\newtheorem{lem}[equation]{Lemma}
\newtheorem{thm}[equation]{Theorem}

\newtheorem{prop}[equation]{Proposition}
\newtheorem{cor}[equation]{Corollary}

\newtheorem{observation}[equation]{Observation}

\theoremstyle{definition}

\theoremstyle{remark}
\newtheorem{remark}[equation]{Remark}

\usepackage{verbatim}
\usepackage{fullpage}
\usepackage[hidelinks]{hyperref} 
\hypersetup{
    colorlinks=true,
}

\renewcommand{\bar}[1]{#1\llap{$\overline{\phantom{\rm#1}}$}}

\begin{document}
\title{Representation complexities of semi-algebraic graphs}
\author{Thao Do}
\address{Department of Mathematics, Massachusetts Institute of Technology, Cambridge MA 02139}
\email{thaodo@mit.edu}

\maketitle
\begin{abstract}
The representation complexity of a bipartite graph $G=(P,Q)$ is the minimum size $\sum_{i=1}^s (|A_i|+|B_i|)$ over all possible ways to write $G$ as a (not necessarily disjoint) union of complete bipartite subgraphs $G=\cup_{i=1}^s A_i\times B_i$ where $A_i\subset P, B_i\subset Q$ for $i=1,\dots, s$. 
 In this paper we prove that if $G$ is \emph{semi-algebraic}, i.e. when $P$ is a set of $m$ points in $\R^{d_1}$, $Q$ is a set of $n$ points in $\R^{d_2}$ and the edges are defined by some semi-algebraic relations, the representation complexity of $G$ is $O( m^{\frac{d_1d_2-d_2}{d_1d_2-1}+\varepsilon} n^{\frac{d_1d_2-d_1}{d_1d_2-1}+\varepsilon}+m^{1+\varepsilon}+n^{1+\varepsilon})$ for arbitrarily small positive $\varepsilon$. This generalizes results by Apfelbaum-Sharir \cite{Ap-Sharir} and Solomon-Sharir \cite{Solomon1,Solomon2,Solomon3}. As a consequence,  when $G$ is $K_{u,u}$-free for some positive integer $u$, its number of edges is $O(u m^{\frac{d_1d_2-d_2}{d_1d_2-1}+\varepsilon} n^{\frac{d_1d_2-d_1}{d_1d_2-1}+\varepsilon}+ u m^{1+\varepsilon}+u n^{1+\varepsilon})$. This bound is stronger than that of Fox, Pach, Sheffer, Suk and Zahl in \cite{Fox} when the first term dominates and $u$ grows with $m,n$. Another consequence is that we can find a large complete bipartite subgraph in a semi-algebraic graph when the number of edges is large.
Similar results hold for semi-algebraic hypergraphs. 
\end{abstract}
\section{Introduction}

Given a bipartite graph $G=(P,Q,\mathcal{E})$ where $P, Q$ are parts and $\mathcal{E}\subset P\times Q$ is the set of edges, we can write $G$ as a union of complete bipartite subgraphs $\E=\cup_{i=1}^s A_i\times B_i$ where $A_i\subset P, B_i\subset Q$ for all $i=1,\dots, s$. The complexity of such a representation of $G$ is defined to be the total number of vertices in the representation $\sum_{i=1}^s (|A_i|+|B_i|)$. The smallest complexity of such a representation, or the \emph{representation complexity} of $G$, is denoted by $J(G)$ or $J(P,Q)$. There is another variation where we require the union to be disjoint, but we will not consider it here.

Compact or compressed representations of graphs have been studied in the algorithmic context before \cite{Visibility graph,graph compression,algorithm application,algorithm application 2}. Representation complexity is similar to the graph's clique covering number (the smallest size of a covering of the graph using cliques or complete bipartite subgraphs), which is shown to be  $\Theta(n^2/\log n)$ for generic graphs on $n$ vertices (see \cite{Chung-E-Spencer,Tuza}). However, for some geometrically defined graphs, it is known that better bounds are possible as there is some structure involved. 
For example, Brass-Knauer \cite{Brass-Knauer} and Apfelbaum-Sharir \cite{Ap-Sharir} found a better upper bound  for the representation complexity  of \emph{point-hyperplane incidence graphs}. Given $m$ points and $n$ hyperplanes in $\R^d$,  their incidence graph is a bipartite graph where one part consists of all the points, the other consists of all the hyperplanes and there is an edge between a point $p$ 
and a hyperplane $H$ if  $p\in H$. 

\begin{thm}[Apfelbaum-Sharir \cite{Ap-Sharir}]\label{Ap-Sharir} The representation complexity of any incidence graph between $m$ points and $n$ hyperplanes in $\R^d$ is $O((mn)^{1-1/d}+m+n)$. 
\end{thm}
In \cite{Solomon1,Solomon2,Solomon3}, Solomon-Sharir found upper bounds for other incidence graphs between points and planes, spheres and other surfaces in $\R^3$. We only state their most general result: a family of curves $\C$ in $\R^3$ is said to have $k$ degrees of freedom
with multiplicity $\mu$, where $k$ and $\mu$ are constants, if (i) for every tuple of $k$ points in $\R^3$ there are at most $\mu$ curves of $\C$ that are incident to all $k$ points and (ii) every pair of curves of $\C$  intersect in at most $\mu$ points. A family $F$ of surfaces in $\R^3$ is said to have $k$ degree of freedom with multiplicity $\mu$ with respect to a given surface $V$
if the family of the irreducible components of the curves \{$\sigma\cap V: \sigma\in F\}$, counted
without
multiplicity, has $k$ degrees of freedom with multiplicity $\mu$ as we just defined.
Theorem 1.7 in \cite{Solomon3} implies the following:
\begin{thm}\label{Solomon-Sharir} The representation complexity of an incidence graph between $m$ points in a bounded degree surface $V$ in $\R^3$ and $n$ algebraic surfaces of bounded degree with $k$ degrees of freedom with respect to $V$ is
$O\left( m^{\frac{k}{2k-1}}n^{\frac{2k-2}{2k-1}}+m+n \right).$
\end{thm}
\subsection{Our result}
In this paper we find an upper bound on the representation complexity of all \emph{semi-algebraic graphs}, which, in some sense, generalizes the above results.


We first define semi-algebraic graphs. Let $G=(P,Q,\mathcal{E})$ be a bipartite graph where $P$ is a set of $n$ points in $\R^{d_1}$ and $Q$ is a set of $m$ points in $\R^{d_2}$. We say $G$ is \emph{semi-algebraic} with description complexity $t$ if there are $t$ polynomials $f_1,\dots, f_t\in \R[x_1,\dots, x_{d_1+d_2}]$, each of degree at most $t$ and a Boolean function $\Phi: \{\pm 1\}^t\to \{\pm 1\}$ such that for any $p\in P, q\in Q$: 
$$ (p,q)\in \mathcal{E} \iff \Phi(f_1(p,q)\geq 0,\dots, f_t(p,q)\geq 0)=1.$$
In other words, we can describe the incidence relation by at most $t$ inequalities involving polynomials of degree at most $t$. We can also view $G$ as an incidence graph between points and semi-algebraic sets. Indeed, for each $p\in P$ and $q\in Q$, define the neighbor sets $$\gamma_p:=\{y\in \R^{d_2}: \phi(f_1(p,y)\geq 0,\dots, f_t(p,y)\geq 0)=1\}$$ 
$$\gamma_q:=\{x\in \R^{d_1}: \phi(f_1(x,q)\geq 0,\dots, f_t(x,q)\geq 0)=1\}.$$
Each $\gamma_p$ is a semi-algebraic set in $\R^{d_2}$, and $\gamma_q$ is a semi-algebraic set in $\R^{d_1}$, both with complexity bounded by $t$. Notice that $(p,q)$ is an edge of $G$ if and only if $q\in\gamma_p$, or equivalently $p\in\gamma_q$. 
Therefore we can view $G$ as the incidence graph between $m$ points $P$ and $n$ semi-algebraic sets $Q^*:=\{\gamma_q:q\in Q\}$ in $\R^{d_1}$, or as the incidence graph between $n$ points $Q$ and $m$ semi-algebraic sets $P^*:=\{\gamma_p:p\in P\}$ in $\R^{d_2}$. When $P^*$ or $Q^*$ consist of algebraic sets (i.e. we can describe the relation by equalities instead of inequalities), we  say $G$ is an \emph{algebraic graph}. The incidence graphs we encounter in Theorems \ref{Ap-Sharir} and \ref{Solomon-Sharir} are both algebraic, and hence also semi-algebraic.

To state our main result, we use the notation $f=O_{a_1,\dots, a_k}(g)$ (or equivalently $f\lesssim_{a_1,\dots,a_k} g$) to denote there is some constant $C$ that depends on $a_1\dots, a_k$ (which we sometimes write $C(a_1,\dots, a_k)$ or $C_{a_1,\dots, a_k})$ such that $f\leq Cg$. 

\begin{thm}\label{main thm}(Main result) Given a bipartite semi-algebraic graph $G=(P,Q,\mathcal{E})$ with description complexity $t$ as above, then for any $\varepsilon>0$,
	\begin{equation}\label{main_ineq}
    J(G)=O_{d_1,d_2,t, \varepsilon} \left(m^{\frac{d_1d_2-d_2}{d_1d_2-1}+\varepsilon} n^{\frac{d_1d_2-d_1}{d_1d_2-1}+\varepsilon}+m^{1+\varepsilon}+n^{1+\varepsilon}\right).
    \end{equation}
    Moreover, if $P$ and $Q$  belong to varieties $V_1\subset \R^{d_1}$ and $V_2\subset \R^{d_2}$ of  degrees at most $t$ and dimensions $e_1$ and $e_2$ respectively, then we can replace $d_i$ by $e_i$ in the above bound, i.e. 
    \begin{equation} \label{main_ineq_2}
    J(G)=O_{d_1,d_2,e_1, e_2, t, \varepsilon} \left(m^{\frac{e_1e_2-e_2}{e_1e_2-1}+\varepsilon} n^{\frac{e_1e_2-e_1}{e_1e_2-1}+\varepsilon}+m^{1+\varepsilon}+n^{1+\varepsilon}\right).
    \end{equation}
    \end{thm}
\begin{remark}
\begin{itemize}
\item  As far as the author is aware, this is the first result on representation complexity of semi-algebraic sets; previous results in \cite{Ap-Sharir,Brass-Knauer, Solomon1,Solomon2,Solomon3} only apply to algebraic sets. 
\item 
   Letting $P$ and $Q$ be sets of points and hyperplanes in $\R^d$ (so $d_1=d_2=d$), we recover Theorem \ref{Ap-Sharir} within an $\varepsilon$ term.  
   \item To recover Theorem \ref{Solomon-Sharir}, let $P$ be the set of points in $V$ (so $e_1=2$) and $Q$ be the set of surfaces in $\R^3$. If all surfaces of $Q$ have degrees at most $t$, $Q$ live in $\R^{{t+3 \choose 3}}$, but since these surfaces have $k$ degrees of freedom with respect to $V$, we expect $e_2=k$. The relationship between the degree of freedom and the dimension of the moduli space of a family of curves/surfaces is quite complicated, see \cite{Solomon3} and the appendix of \cite{Sharir-Zahl} for more details. But typically we should expect $e_2=k$  and hence recover the bound in Theorem \ref{Solomon-Sharir} within an $\varepsilon$ term.
   \item A similar result holds for  \emph{semi-algebraic hypergraphs}, see Theorem \ref{rep comp hypergraph} for details.
 \end{itemize}
 \end{remark}

 \subsection{Applications}
We first present an application to the Zarankiewicz's problem \cite{Zarankiewcz}, a central problem in graph theory and incidence geometry. It  asks for the largest possible number of edges in an $m\times n$  bipartite graph that avoids $K_{u,u}$ for some fixed positive integer $u$. Here $K_{u,u}$ denotes the complete bipartite graph of size $u\times u$, and we say a graph $G$ avoids $H$ or $G$ is $H$-free if $G$ does not contain any subgraph congruent to $H$.  In 1954, K\H{o}v\'ari, S\'os and Tur\'an proved a general upper bound of form $O_u(mn^{1-1/u}+n)$, which is only known to be tight for $u=2$ and $u=3$. 
 
 Better bounds are known when the graphs are incidence graphs between points and lines in $\R^2$ (Szemer\'edi-Trotter Theorem \cite{S-T}), points and curves in $\R^2$ (Pach-Sharir \cite{Pach-Sharir}), and points and hyperplanes in $\R^d$ (Apfelbaum- Sharir \cite{Ap-Sharir}). Recently, Fox, Pach, Sheffer, Suk and Zahl \cite{Fox} generalized these results to all semi-algebraic graphs.

\begin{thm}[Fox, Pach, Sheffer, Suk and Zahl  \cite{Fox}]\label{Fox}
Given a bipartite semi-algebraic graph $G=(P,Q,\mathcal{E})$ with description complexity $t$ as above, if $G$ avoids $K_{u,u}$ then for any $\varepsilon>0$, $$|\mathcal{E}(G)|= O_{t,d_1,d_2,u,\varepsilon}\left(m^{\frac{d_1d_2-d_2}{d_1d_2-1}+\varepsilon}n^{\frac{d_1d_2-d_1}{d_1d_2-1}}+m+n\right).$$ 
When $d_1=d_2=2$ we can delete the $\varepsilon$ term. 
\end{thm}

This theorem assumes $u$ is a fixed constant and does not explicitly state how the bound depends on $u$. Following its proof in \cite{Fox}, we estimate\footnote{This is done by keeping track of the dependence of $u$ in each step of the proof. It is possible that this bound can be improved using a more careful analysis, but to the author it seems quite infeasible to obtain anything where the dependence of $u$ in the first term is better than linear. } the bound to be 
\begin{equation}\label{Fox et.al. bound}_{d_1,d_2,\varepsilon, t} (u^{1+\frac{(d_1-1)(d_2+1)}{d_1d_2-1}}m^{\frac{d_1d_2-d_2}{d_1d_2-1}+\varepsilon} n^{\frac{d_1d_2-d_1}{d_1d_2-1}}+um+un).\end{equation}
Using our Theorem \ref{main thm} we get the following bound.
\begin{cor}\label{cor-Z-problem} Assume the given semi-algebraic graph $G$ is $K_{u,u}$-free, then
\begin{equation}\label{Z prob our bound}
|\mathcal{E}(G)|=O_{d_1,d_2,t,\varepsilon}\left( um^{\frac{d_1d_2-d_2}{d_1d_2-1}+\varepsilon} n^{\frac{d_1d_2-d_1}{d_1d_2-1}+\varepsilon}+u m^{1+\varepsilon}+u n^{1+\varepsilon}\right).
\end{equation}
\end{cor}
 This follows from the fact that $|\E|\leq (s+t) J(G)$ whenever the graph $G$ is $K_{s,t}$-free. Indeed, in any decomposition $\E=\cup A_i\times B_i$,  for each $i$ either $|A_i|<s$ or $|B_i|<t$. In either case $|A_i||B_i|\leq (s+t)(|A_i|+|B_i|)$, taking the sum over all $i$ we get $|\E|\leq (s+t) J(G)$.
 
In \eqref{Fox et.al. bound}, when  $m>n^{d_1}$ or $n>m^{d_2}$, the dominant term in $m+n$, in which our bound \eqref{Z prob our bound} is $\varepsilon$ weaker. On the other hand, when $n^{1/d_2}\leq m\leq n^{d_1}$, which is usually the most interesting range, the first term dominates and when $u$ gets large, the bound in \eqref{Z prob our bound} is stronger.

Can we prove a sub-linear dependence on $u$ in the first term of \eqref{Z prob our bound}? 
If  the graph is $K_{s,u}$-free where $s$ is fixed and $u$ can get large, the answer is yes.

\begin{thm}\label{no K_2,u}[Lund, Sheffer and de Zeeuw \cite{bisector energy}]
Given $m$ points and $n$ varieties of degree at most $t$ in $\R^d$ such that their incidence graph is $K_{s,u}$-free where $s$ is a fixed small integer and $u$ can be large. Then the number of point-variety incidences is at most 
$$O_{d,t,s,\varepsilon}(m^{\frac{ds-s}{ds-1}+\varepsilon} n^{\frac{ds-d}{ds-1}}u^{\frac{d-1}{ds-1}}+n+um).$$

\end{thm}
 
 However, in this case  the condition that $s$ is fixed is crucial because $s$ appears in the exponent of $m,n$. Therefore this result is not as robust as Corollary \ref{cor-Z-problem}. 

On the other hand, here is an heuristic argument why we should not expect to improve the dependence of $u$ in the first term in \eqref{Z prob our bound} better than $u^{2/3}$. Szemer\'edi-Trotter's theorem \cite{S-T} is known to be tight: there exist $m$ points and $n$ lines in the plane with $\Theta(m^{2/3}n^{2/3}+m+n)$ incidences. If we replace each point and each line by $u-1$ copies of them, we have a configuration with $(u-1)m$ points, $(u-1)n$ lines, no $K_{u,u}$ and $\Theta((u-1)^2 m^{2/3}n^{2/3}+(u-1)^2m+(u-1)^2 n)=\Theta (u^{2/3} (um)^{2/3}(un)^{2/3}+u(um)+u(un))$ incidences.

The next application is to find a large complete bipartite subgraph in an incidence graph when the number of edges is large.


In \cite{Ap-Sharir}, Apfelbaum and Sharir explored the following question: when $m$ points and $n$ hyperplanes in $\R^d$ form many incidences, what can we say about the size of the largest complete subgraph found in their incidence graph? 
This question is related to Ramsey theory. It is known that for any graph $G$ on $n$ vertices, there are $P,Q$ of size $\log n$ such that $P\times Q$ is fully contained in $G$ or its complement graph $\bar{G}$. Stronger results hold for semi-algebraic graphs: if a semi-algebraic graph $G$ with description complexity $t$ on $n$ vertices in $\R^d$ has $\Theta(n^2)$ edges, we can find a complete bipartite subgraph of size $\Theta(n)\times \Theta(n)$  (see  \cite{Semi-hypergraph regularity lemma, Semi-hypergraph regularity lemma2}).
However, not much is known when the number of edges is neither too small nor too large.

Using our main theorem, we easily obtain the following result:
  \begin{cor}
 If $I=\Omega( m^{\frac{d_1d_2-d_2}{d_1d_2-1}+\varepsilon} n^{\frac{d_1d_2-d_1}{d_1d_2-1}+\varepsilon}+m^{1+\varepsilon}+n^{1+\varepsilon})$ then $G$ contains some $K_{u_1,u_2}$ with 
 $$u_1u_2\geq \min  \left\{\left(\frac{I}{ m^{\frac{d_1d_2-d_2}{d_1d_2-1}+\varepsilon}n^{\frac{d_1d_2-d_1}{d_1d_2-1}+\varepsilon}}\right)^2,\frac{I^2}{n^{2+2\varepsilon}}, \frac{I^2}{m^{2+2\varepsilon}},\right\}.$$
 \end{cor}
 Indeed, consider the most compact decomposition $\E=\cup A_i\times B_i$. Since $J(G)=\sum (|A_i|+|B_i|)$ and $I(G)=\sum |A_i||B_i|$, there exists some $i$ such that $ \frac{|A_i||B_i|}{A_i+B_i}\geq \frac{I}{J(G)}$. Let $u_1=|A_1|, u_2=|B_1|$, then each $u_i$ is at least $I/J(G)$, and using  \eqref{main_ineq} we get our desired result.

\subsection{Organization}

We first recall two main tools used in the proof in section \ref{sec:preliminary}: a Milnor-Thom type result and the polynomial partitioning method. The main theorem is proved in section \ref{sec:proof}.
In section \ref{sec:hypergraph}, we extend the result to \emph{semi-algebraic hypergraphs}. 
Finally we end with some open questions in section \ref{sec:conclusion}.
 \subsection{Acknowledgement} The author would like to thank Larry Guth for his support and guidance throughout the project. She also thanks Esther Ezra for mentioning \cite{range searching} which inspired the proof of Proposition \ref{step 1} and thanks Adam Sheffer for pointing out the reference of Theorem \ref{no K_2,u}. Finally, she thanks Malcah Effron and Vishesh Jain for proofreading the preprint.

 \section{Preliminary}\label{sec:preliminary}
 \subsection{Milnor-Thom type results}
  Milnor-Thom's theorem \cite{Milnor, Thom} states that the zero set  of a degree $D$ polynomial $f$, denoted by $Z(f)$, divides $\R^d$  into at most $(50D)^d$ connected components (i.e. $\R^d\setminus Z(f)$ has at most $(50D)^d$ connected components). Basu, Pollack and Roy extended this result to the case when we restrict our attention to a variety inside $\R^d$. 
 
 A \emph{sign pattern} for a set of $s$ $d$-variate polynomials $\{f_1,\dots, f_s\}$ is a vector $\sigma\in \{-1,0,+1\}^s$. A sign pattern $\sigma$ is \emph{realizable} over a variety $V\subset\R^d$ if there is some $x\in V$ such that (sign$(f_1(x))$, sign$(f_2(x)),\dots,$sign$(f_s(x)))=\sigma$. The set of all such $x$ is the realization space of $\sigma$ in $V$, denoted by $\Omega_\sigma$.
\begin{thm}[Basu, Pollack and Roy, 1996 \cite{Basu-Pollack-Roy}]\label{Basu-Pollack-Roy}
Given positive integers $d, e, M, t, l$, let $V$ be an $e$-dimensional real algebraic set in $\R^d$ of complexity\footnote{A variety has \emph{complexity} at most $M$ if it can be realized as the intersection of zero-sets of at most $M$ polynomials, each of degree at most $M$.} at most $M$, 
and let $f_1,\dots, f_s$ be $d$-variate real polynomials of degree at most $t$. Then the total number of connected components of $\Omega_\sigma$ for all realizable sign patterns $\sigma$ of $\{f_1,\dots, f_s\}$ is at most $O_{M,d,e}((ts)^e)$. 
\end{thm}

This result  implies if we restrict to a variety $V$ with dimension $e$ and bounded complexity in $\R^d$, then the number of connected components that $f_1,\dots, f_s$ partition $V$ grows with $e$ instead of $d$. 
Furthermore, a similar result holds if we replace $V$ by $V\setminus W$ for some variety $W$ with bounded complexity.
\begin{thm}[Theorem A.2 in \cite{Tao}] 
\label{Tao}
Given positive integers $d, e, M, t$ such that $e\leq d$, let $V$ and $W$ be a real algebraic sets in $\R^d$ of complexity at most $M$ such that $V$ is $e$-dimensional. Then for any polynomial $P:\R^d\to \R$ of degree $t \geq 1$, the set $\{x\in V\setminus W: P(x)\neq 0\}$ has $O_{M,d,e}(t^e)$ connected components.
\end{thm}
\subsection{Polynomial partitioning}\label{sec:poly partitioning}
Polynomial partitioning method was first introduced by Guth and Katz in \cite{Larry} in 2010 and developed to several different versions since then. In this paper we use the version proved in \cite{Fox}. 
Given $n$ points in $\R^d$, we say a polynomial $f\in\R[x_1,\dots,x_d]$ is an $r$-partitioning if the zero set of $f$, denoted by $Z(f)$, divides the space into open connected components and each component contains at most $n/r$ points of the given.

\begin{thm}[Theorem 4.2 in \cite{Fox}]\label{poly partitioning}
Let $P$ be a set of points in $\R^d$, and let $V\subset \R^d$ be an irreducible variety of degree $D$ and dimension $d'$. Then for big enough $r$, there exists an $r-$partitioning polynomial $g$ for $P$ such that $g\notin I(V)$ and $\deg g\leq C_{part} \cdot r^{1/d'}$ where $C_{part}$ depends only on $d$ and $D$. 
\end{thm}
This theorem implies in $\R^d$, if we restrict our attention to points in an irreducible variety of small degree and dimension $d'<d$, then we can perform a polynomial partitioning the same way as in $\R^{d'}$.


\section{Proof of the main theorem}\label{sec:proof}

We follow the general strategy used to prove Theorem \ref{Fox} and Theorem \ref{no K_2,u}: first obtain a nontrivial bound by various ways, then use polynomial partitioning to get the desired stronger bound.  
In particular, in Theorem \ref{Fox}, Fox, Pach, Sheffer, Suk and Zahl first proved $|\E|= O(mn^{1-1/d_2}+n)$ using a packing result in VC-dim theory, and then used polynomial partitioning in the space $\R^{d_1}$ (where the point set $P$ lives in). In Theorem \ref{no K_2,u}, Lund, Shefer and de Zeeuw first showed that $|\E|=O(mn^{1-1/s}+n)$ by K\H{o}v\'ari, S\'os and Tur\'an, then applied polynomial partitioning in $\R^d$. Notice the similar roles between $s$ and $d_2$ in these two theorems, which is partially explained in Remark \ref{remark after main thm}.  The main difficulty for us is that we do not have the condition of $K_{u,u}$-free which is crucial in K\H{o}v\'ari, S\'os and Tur\'an's inequality and the VC-dim argument. We fix this by using polynomial partitioning in $\R^{d_2}$ (the space where the point set $Q$ lives in) to get a slightly weaker bound $O(mn^{1-1/d_2+\varepsilon}+n^{1+\varepsilon})$. This is why in our final result we have the term $n^{1+\varepsilon}+m^{1+\varepsilon}$ instead of $m+n$. Another main difference with previous proofs is one case in the analysis of incidences (or edges) after applying polynomial partitioning: incidences between points in a cell (of the partitioning) and  semi-algebraic sets that contain that cell. In  Theorem \ref{Fox} and \ref{no K_2,u}, since there is no $K_{u,u}$ or $K_{s,u}$, the total number of edges in this case is small. In our situation, the number of edges can be large but the representation complexity remains small since these edges form complete bipartite subgraphs.

We start with the first step.

\begin{prop}\label{step 1} Given a semi-algebraic bipartite graph $G=(P,Q,\E)$ with description complexity $t$ where $P$ is a set of $m$ points in $\R^{d_1}$ and $Q$ is a set of $n$ points in $\R^{d_2}$ then $J(G)= O_{d_1,d_2,t,\varepsilon} (mn^{1-1/d_2+\varepsilon}+n^{1+\varepsilon})$ for arbitrarily small $\varepsilon$.
\end{prop}

\begin{remark} This is somewhat related to a result by Agarwal, Matousek and Sharir on range searching with semi-algebraic sets \cite{range searching}. Given $n$ points in $\R^d$, the range searching problem asks for a way to determine how many points a semi-algebraic set contains.  They proved this problem can be solved with $O(n)$ storage,  $O(n\log n)$ expected processing time and  $O(n^{1-1/d+\varepsilon})$ query time. This means if there are $m$ semi-algebraic set, we need $O(mn^{1-1/d+\varepsilon}+n\log n)$ time to learn about their point-set incidence structure. This is almost the same with our bound of representation complexity of their incidence graph. In fact with some extra work we can prove $J(G)=O(mn^{1-1/d_2+\varepsilon}+n\log n)$, but we do not include a proof here because it is not necessary for our problem. 
\end{remark}

\begin{proof} We shall prove a more general statement: If $Q$ belongs to an irreducible variety $V$ of dimension $e_2$ (for some $e_2\leq d_2$) and degree $D$, then there exists some constant $C^*_{e_2}$ that depends on $d_1,d_2,e_2,t,D$ and $\varepsilon$ such that
$$J(G)\leq C^*_{e_2} (mn^{1-1/e_2+\varepsilon}+n^{1+\varepsilon}).$$
We prove this by induction on $e_2$ and $m+n$. The statement is vacuous when $e_2=0$. When $m+n$ is small, we can choose the constant big enough for the inequality to hold true. 


For the induction step, we use polynomial partitioning with respect to the variety $V$. Let $r$ be a parameter to be chosen later. By Theorem \ref{poly partitioning}, there exists a polynomial $f$ of degree at most $r$ that partition $V$ into $s\leq c_1 r^{e_2}$ cells $\Omega_1,\dots, \Omega_s$, each contains at most $ c_2 n/r^{e_2}$ points of $Q$. Here $c_1,c_2$ are constants that depends on $d_1,d_2, e_2$ and $D$. 
For each $i$, let $Q_i, P_i, P_i'$ respectively denote the set of points of $Q$ contained in $\Omega_i$, the set of semi-algebraic sets in $P^*$ that contains $\Omega_i$ and the set of semi-algebraic sets in $P^*$ that \emph{crosses} $\Omega_i$ (i.e. have nonempty intersection but not contain). Finally let $Q'$ denote the set of points of $Q$ that lies on $Z(f)$. More precisely:
$$Q_i=Q\cap \Omega_i;\quad Q':=Q\cap Z(f);$$
$$P_i:=\{\gamma\in P^*: \Omega_i\subset\gamma\};\quad
P_i':=\{\gamma\in P^*:\gamma\cap \Omega_i\neq\emptyset, \Omega_i\not\subset \gamma\}.$$

We can decompose our graph as 
$G(P^*, Q)=\cup_{i=1}^s G(P_i, Q_i)\cup_{i=1}^s G(P'_i, Q_i)\cup G(P^*,Q')$, which implies
 \begin{equation} \label{decompose J}
J(G)\leq \sum_{i=1}^s J(P_i, Q_i)+\sum_{i=1}^s J(P_i',Q_i)+ J(P^*, Q').
\end{equation}
We shall bound each term on the RHS of \eqref{decompose J} by $1/3C^*_{e_2} (mn^{1-1/e_2+\varepsilon}+n^{1+\varepsilon})$ for appropriate choices of $r$ and $C^*_{e_2}$, which would complete our induction step.
To do so, we first make some observations:
\begin{observation}\label{sum p_i q_i} \begin{itemize}
\item [(i)] $|Q'|+\sum_{i=1}^s |Q_i|=n$ and
$|Q_i|\leq c_2 n/r^{e_2}$ for each $i$.
\item[(ii)] 
 $\sum_{i=1}^s |P_i|= \#\{(\Omega, \gamma):\Omega\subset \gamma\} \leq s m\leq c_1 r^{e_2} m$
 \item [(iii)] $\sum_{i=1}^s |P_i'|=\#\{(\Omega, \gamma):\gamma\text{ crosses } \Omega\} \leq c_3 mr^{e_2-1}$ for some constant $c_3$ that depends on $D,d_1,e_2,t$. 
\end{itemize}
\end{observation}
Only the last inequality requires some reasoning: 
We claim that each semi-algebraic set $\gamma_p$ crosses  at most $O(r^{e_1-1})$ cells. Indeed, each $\gamma_p$ is defined by $t$ polynomials $f_1(p, x),\dots, f_t(p,x)$. In order for $\gamma_p$ to cross a cell in $V\setminus Z(f)$, some polynomial, say $f_1$, must not vanish on $V$. Then $Z(f_1)\cap V$ is some variety of dimension at most $e_1-1$. By theorem \ref{Tao}, $f$ partitions this variety in at most $O_{t,d_2, D} (r^{e_1-1})$ cells; this in turn implies $\gamma_p$ crosses $O(r^{e_2-1})$ cells. Adding them up, we get the inequality in (iii). 
\\
\\
\noindent\textbf{Bounding the first term in \eqref{decompose J}}: Since each set in $P_i$ contains the cell $\Omega_i$ which contains $Q_i$, $G(P_i,Q_i)$ is a complete graph. Thus $\sum_{i=1}^s J(P_i,Q_i)= \sum_{i=1}^s (|P_i|+|Q_i|)\leq c_2(r^{e_2}m+n)$ by observation \ref{sum p_i q_i}. This is bounded by $1/3 C^*_{e_2}( mn^{1-1/d_2+\varepsilon}+n^{1+\varepsilon})$ if we  choose $C_{e_2}^*>3c_1$ and $r^{e_2}<n^{1-1/d_2+\varepsilon}$.
\\
\\
\textbf{Bounding the second term:} Since each cell contains fewer than $n$ points, we can apply the induction assumption to each cell, sum them up and use observation \ref{sum p_i q_i}:

\begin{align*}
\sum_{i=1}^s J(P_i',Q_i)& \leq \sum C^*_{e_2}(|P_i'||Q_i|^{1-1/e_2+\varepsilon}+|Q_i|^{1+\varepsilon})\nonumber\\
&\leq C_{e_2}^*\sum_{i=1}^s \left[|P_i'| \left(\frac{c_2n}{r^{e_2}}\right)^{1-1/e_2+\varepsilon}+\left(\frac{c_2n}{r^{e_2}}\right)^{1+\varepsilon}\right]\nonumber\\
&\leq C_{e_2}^*\left(m r^{e_2-1}\left(\frac{c_2n}{r^{e_2}}\right)^{1-1/e_2+\varepsilon} + c_1r^{e_2}\left(\frac{c_2n}{r^{e_2}}\right)^{1+\varepsilon}\right)\nonumber\\
& \leq C_{e_2}^* r^{-e_2\varepsilon}\left(c_2^{1-1/e_2+\varepsilon} mn^{1-1/e_2+\varepsilon}  + c_1c_2^{1+\varepsilon} n^{1+\varepsilon} \right)
\end{align*}
By choosing $r$ big enough compared to $c_1,c_2$ we get our desired bound.
\\
\\
\noindent \textbf{Bounding the third term:} We notice that $Q'$ belong to $V\cap Z(f)$. Since $V$ is irreducible of dimension $e_2$ and degree $D$, its intersection with $Z(f)$ must be disjoint union of at most $Dr$ irreducible varieties of dimension at most $e_2-1$ and degree bounded by $D$ and $r$ (since $\deg f\leq r)$. Apply our induction assumption to each of those irreducible varieties, we get: $J(P^*, Q')\leq Dr C^*_{e_2-1} (mn^{1-1/(e_2-1)+\varepsilon}+n^{1+\varepsilon})\leq 1/3 C^*_{e_2}(mn^{1-1/e_2+\varepsilon}+ n^{1+\varepsilon})$ by choosing $C^*_{e_2}>3DC^*_{e_2-1}$ and $r<n^{1/e_2(e_2-1)}$.

This completes our proof of Proposition \ref{step 1}.
\end{proof}

For the second step, we prove \eqref{main_ineq_2} (the stronger statement involving $e_i$) by using polynomial partitioning again but in a different space $\R^{d_1}$, in which we view $G$ as the incidence graph between $P$ and $Q^*$.

We proceed similarly to the proof of Proposition \ref{step 1}. To simplify in many cases we write $f\lesssim g$ or $f=O(g)$ without explicitly stating what the constant depends on. 
Fix $e_2$. By Proposition \ref{step 1} we know $J(G)=O(mn^{1-1/e_2+\varepsilon}+n^{1+\varepsilon})$. If $n\geq m^{e_2}$, $n^{1+\varepsilon}$ is the main term in the right hand side of Proposition \ref{step 1}, which implies $J(G)\lesssim n^{1+\varepsilon}$ and hence \eqref{main_ineq_2} holds. Similarly by symmetry $J(G)=O(nm^{1-1/e_1+\varepsilon}+m^{1+\varepsilon})$ and \eqref{main_ineq_2} holds when $m\geq n^{e_1}$.
Hence from now on we can assume $n\leq m^{e_2}$ and $m\leq n^{e_1}$. In this case the term $m^{\frac{e_1e_2-e_2}{e_1e_2-1}+\varepsilon}n^{\frac{e_1e_2-e_1}{e_1e_2-1}+\varepsilon}$ is the dominant term in our bound.  

We use induction by $e_1$ and $m+n$. The statement is vacuous when $e_1=0$. When $m+n$ is small, we can choose the constant big enough for the inequality to hold true. 
For the induction step, we use polynomial partitioning with respect to the variety $V_1$. Let $r$ be a parameter to be chosen later. By theorem \ref{poly partitioning}, there exists a polynomial $f$ of degree at most $r$ that partition $V_1$ into $s=O(r^e_1)$ cells $\Omega_1,\dots, \Omega_s$, each contains $O(n/r^{e_1})$ points of $P$.

For each $i$, let $P_i, Q_i, Q_i'$ respectively denote the set of points of $P$ contained in $\Omega_i$, the set of semi-algebraic sets in $Q^*$ that contains $\Omega_i$ and the set of semi-algebraic sets in $Q^*$ that crosses $\Omega_i$ (i.e. have nonempty intersection but not contain). Finally let $P'$ denote the set of points of $P$ that lies on $Z(f)$. Similar to Observation \ref{sum p_i q_i} we have:
\begin{itemize}
\item [(i)] $\sum_{i=1}^s |P_i|+|P'|=m$ and $|P_i|\lesssim n/r^{e_1}$ for each $i$.
\item[(ii)] 
 $\sum_{i=1}^s |Q_i|\lesssim  r^{e_1} m$
 \item[(iii)] $\sum_{i=1}^s |Q_i'|\lesssim mr^{e_1-1}$.
\end{itemize}

We can decompose our graph as 
$G(P^*, Q)=\cup_{i=1}^s G(P_i, Q_i)\cup_{i=1}^s G(P_i, Q'_i)\cup G(P',Q^*)$, which implies
 \begin{equation}\label{decompose J2}
J(G)\leq \sum_{i=1}^s J(P_i, Q_i)+\sum_{i=1}^s J(P_i,Q_i')+ J(P', Q^*).
\end{equation}
We now bound each term in the RHS of \ref{decompose J2} by $C(m^{\frac{e_1e_2-e_2}{e_1e_2-1}+\varepsilon}n^{\frac{e_1e_2-e_1}{e_1e_2-1}+\varepsilon}+m^{1+\varepsilon}+n^{1+\varepsilon}$ for appropriate choice of $C$ and $r$.
The first term is bounded by the exact same way as before: since $G(P_i,Q_i)$ is a complete graph. Thus $\sum_{i=1}^s J(P_i,Q_i)= \sum_{i=1}^s (|P_i|+|Q_i|)= O(r^{e_1}n+m)=O(n^{1+\varepsilon}+m^{1+\varepsilon})$  by choosing $r$ so that $r^{e_1}\lesssim n^\varepsilon$. 

For the second term: we apply the induction assumption in each cell, sum them up, and then use H\"older's inequality
\begin{align*}
\sum_{i=1}^s J(P_i,Q'_i)& \lesssim \sum (|P_i|^{\frac{e_1e_2-e_2}{e_1e_2-1}+\varepsilon}|Q'_i|^{\frac{e_1e_2-e_1}{e_1e_2-1}+\varepsilon}+|P_i|^{1+\varepsilon}+|Q'_i|^{1+\varepsilon})\nonumber\\
&\lesssim \left(\frac{m}{r^{e_1}}\right)^{\frac{e_1e_2-e_2}{e_1e_2-1}+\varepsilon}\sum |Q'_i|^{\frac{e_1e_2-e_1}{e_1e_2-1}+\varepsilon}+s\left(\frac{m}{r^{e_1}}\right)^{1+\varepsilon}+(\sum |Q'_i|)^{1+\varepsilon}\nonumber\\
&\lesssim  r^{e_1} \left(\frac{m}{r^{e_1}}\right)^{\frac{e_1e_2-e_2}{e_1e_2-1}+\varepsilon} \left(\frac{n}{r}\right)^{\frac{e_1e_2-e_1}{e_1e_2-1}+\varepsilon}+ m^{1+\varepsilon}r^{-e_1\varepsilon}+n^{1+\varepsilon}r^{(e_1-1)(1+\varepsilon)}
\nonumber\\
& \lesssim r^{-(e_1+1)\varepsilon} m^{\frac{e_1e_2-e_2}{e_1e_2-1}+\varepsilon}n^{\frac{e_1e_2-e_1}{e_1e_2-1}+\varepsilon}+r^{-e_1\varepsilon} m^{1+\varepsilon}+n^{1+\varepsilon}r^{(e_1-1)(1+\varepsilon)}
\end{align*}
By our assumption $m^{e_1}>n$ at the beginning, $m^{\frac{e_1e_2-e_2}{e_1e_2-1}}n^{\frac{e_1e_2-e_1}{e_1e_2-1}+\varepsilon}\gtrsim n^{1+\varepsilon}$. Hence we obtain the bound by choosing $r$ not too large such that $r^{(e_1-1)(1+\varepsilon)}< r^{-\varepsilon(e_1+1)} m^\varepsilon$ and not too small so that $r^{-(e_1+1)\varepsilon}$ outweights the constant.

For the third term: notice that $P'$ belong to $V_1\cap Z(f)$. Since $V_1$ is irreducible of dimension $e_1$ and degree $D$, its intersection with $Z(f)$ must be disjoint union of several irreducible varieties of dimension at most $e_1-1$ and degree bounded by $D$ and $r$ (since $\deg f\leq r)$. Apply the induction assumption for $e_1-1$: 
$$J(P', Q^*)\lesssim m^{\frac{(e_1-1)e_2-e_2}{(e_1-1)e_2-1}+\varepsilon}n^{\frac{(e_1-1)e_2-(e_1-1)}{(e_1-1)e_2-1}+\varepsilon}+n^{1+\varepsilon}+m^{1+\varepsilon}.$$
By simple algebra, $m^{\frac{(e_1-1)e_2-e_2}{(e_1-1)e_2-1}+\varepsilon}n^{\frac{(e_1-1)e_2-(e_1-1)}{(e_1-1)e_2-1}+\varepsilon}\leq m^{\frac{e_1e_2-e_2}{e_1e_2-1}+\varepsilon}n^{\frac{e_1e_2-e_1}{e_1e_2-1}+\varepsilon}$ when $m<n^{e_1}$ and $n< m^{e_2}$ which holds true by our assumption at the beginning. Hence this term is bounded as we wished. This completes the proof of theorem \ref{main thm}.

\qed

\section{Extension to semi-algebraic hypergraphs}\label{sec:hypergraph}
Theorem \ref{main thm} generalizes naturally to semi-algebraic hypergraphs.
A hypergraph $H$ is called $k-$uniform  if each hyperedge is a $k$-tuple of its vertices. It is $k$-partite if its vertices can be partitioned into $k$ disjoint subset $P_1,\dots, P_k$ and each hyperedge is some tuple $(p_1,\dots, p_k)$ where $p_i\in P_i$ for $i=1,\dots, k$. We usually use $\mathcal{E}$, or $\mathcal{E}(H)$ to denote the set of hyperedges of $H$. 

Let $H$ be a $k$-uniform $k$-partite hypergraph $H=(P_1,\dots, P_k, \mathcal{E})$ where $P_i$ is a set of $n_i$ points in $\R^{d_i}$ for $i=1,\dots, k$ and $\mathcal{E}$ is the set of all hyperedges. This hypergraph is said to be \emph{semi-algebraic with description complexity $t$} if there are $t$ polynomials $f_1,\dots,f_t\in \R[x_1,\dots, x_{d_1+\dots+d_k}]$, each of degree at most $t$, and a Boolean function $\Phi(X_1,\dots, X_t)$ such that for any $p_i\in P_i$, $i=1\dots, k$:
$$(p_1,\dots, p_k)\in \mathcal{E} \iff \Phi(f_1(p_1,\dots, p_k)\geq 0,\dots, f_t(p_1,\dots, p_k)\geq 0)=1.$$

Semi-algebraic hypergraphs have been studied extensively recently (see for example \cite{Ramsey hypergraph, Semi-hypergraph regularity lemma2,Semi-hypergraph regularity lemma}). Many classical results about hypergraphs such as the Ramsey's bound and Szemer\'edi's regularity lemma can be improved in the semi-algebraic setting. 

Recently the author  extends theorem \ref{Fox} to semi-algebraic hypergraph in \cite{Z problem hypergraph}.
To state that result, we need some definitions.
Let $\vec{d}=(d_1,\dots, d_k)$ and $\vec{n}=(n_1,\dots, n_k)$ be  vectors in $\mathbb{Z}^k$. 
For each $\vec{d}$ such that $d_i\geq 1$ for all $i$, and each $\varepsilon>0$, define functions $E_{\vec{d}}(\vec{n})$ and $F^\varepsilon_{\vec{d}}:\R^k\to R$ as followed:
\begin{equation}\label{defi_E}
E_{\vec{d}}(\vec{n})=E_{d_1,\dots, d_k}(n_1,\dots, n_k):= \prod_{i=1}^k n_i^{1-\frac{1/(d_i-1)}{k-1+\frac{1}{d_1-1}+\dots+\frac{1}{d_k-1}}}.
\end{equation}
\begin{equation}\label{defi_F}
F^\varepsilon_{\vec{d}}(\vec{n}):=\sum_{I\subset [k], |I|\geq 2}  E_{\vec{d_I}}(\vec{n_I})\prod_{i\in I} n_i^\varepsilon \prod_{i\notin I} n_i  + \left(\frac{1}{n_1}+\dots+\frac{1}{n_k}\right)\prod_{i=1}^k n_i
\end{equation}
Notice that $E_{d_1,d_2}(m,n)= m^\frac{d_1d_2-d_2}{d_1d_2-1}n^\frac{d_1d_2-d_1}{d_1d_2-1}$ and $F^\varepsilon_{d_1,d_2}(m,n)$ is exactly the bound in theorem \ref{Fox}. For properties of those functions $E,F$, see appendix \ref{property E and F}. 
Here is the main result in \cite{Z problem hypergraph}:
\begin{thm}[Do \cite{Z problem hypergraph}]\label{Z hypergraph}
Given a $k$-uniform $k$-partite hypergraph $H=(P_1,\dots, P_k,\mathcal{E})$ with description complexity $t$ as above, if $H$ avoids $K_{u,\dots,u}$ for some fixed $u$ then 
$$|\mathcal{E}(H)|=O_{t,k,u, \vec{d}, \varepsilon}\left(F^\varepsilon_{\vec{d}}(\vec{n})\right).$$ 
Moreover, if for each $i\leq k$, $P_i$ belongs to an irreducible variety of degree $D$ and dimension $e_i\leq d_i$, then 
$|\mathcal{E}(H)|=O_{t,k,u, \vec{d}, D, \varepsilon}\left(F^\varepsilon_{\vec{e}}(\vec{n})\right)$
where $\vec{e}=(e_1,\dots, e_{k})$.
\end{thm}

Given a $k$-uniform $k$-partite hypergraph $H=(P_1,\dots,P_k,\E)$, we define its representation complexity  as followed. For each way of decomposing it as disjoint union of $s$ complete $k$-partite subhypergraphs $\E=\cup_{i=1}^s A_{i1}\times A_{i2}\times \dots\times A_{ik}$ where $A_{ij}\subset P_j$ for all $i,j$, its complexity is $$\sum_{i=1}^s |A_{i1}||A_{i2}|\dots|A_{ik}|\left(\frac{1}{|A_{i1}|}+\frac{1}{|A_{i2}|}+\dots+\frac{1}{|A_{ik}|}\right).$$ The smallest such quantity among all decompositions is called the \emph{representation complexity} of $H$, denoted by $J(H)$.
\begin{remark}
It might seem natural to define complexity as $\sum_{i=1}^s (|A_{i1}|+\dots+|A_{ik}|)$ since we only need the information about vertices of $A_{ij}$ to represent $H$. However, our definition has the advantage of preserving the desired property: $J(H)\leq k|\E|$ and when $H$ is $K_{u_1,\dots, u_k}$-free, $J(H)\geq (u_1+\dots+u_k) |\mathcal{E}|$. 
\end{remark}

In this paper we shall prove the following result: 
\begin{thm}\label{rep comp hypergraph} [Representation complexity of semi-algebraic hypergraphs]
Given a $k$-uniform $k$-partite hypergraph $H=(P_1,\dots, P_k,\mathcal{E})$ with description complexity $t$ as above. Then 
$$J(H)=O_{t,k,u, \vec{d}, \varepsilon}\left(F^{*,\varepsilon}_{\vec{d}}(\vec{n})\right)$$
where $$F^{*,\varepsilon}_{\vec{d}}(\vec{n}):=\prod_{i=1}^k n_i^{\varepsilon} \sum_{I\subset[k]}\left(E_{\vec{d_I}}(\vec{n_I})\prod_{i\notin I} n_i\right).$$
Moreover, if for each $i\leq k$, $P_i$ belongs to an irreducible variety of degree $D$ and dimension $e_i\leq d_i$ then we can replace $d_i$ by $e_i$ respectively. 
\end{thm}
As a corollary, when the hypergraph is $K_{u,\dots, u}$-free we get a bound on $|\E(H)|$. 
$$|\E(H)|= O_{t,k, \vec{d}, \varepsilon}u\left(F^{*,\varepsilon}_{\vec{d}}(\vec{n})\right)$$
Similar with the graph case, this bound is $\varepsilon$-weaker than that in Theorem \ref{Z hypergraph} for certain range of the $n_i's$, but have a better dependence in $u$  when the first term dominates. 

Since the proof is almost identical with that in the previous section with the only additional new idea of the \emph{grid polynomial partitioning} developed in \cite{Z problem hypergraph}, we shall only give a sketch of the proof.
\begin{proof}[Sketch of the proof] 
We only prove the stronger statement involving $e_i$.
We follow the same strategy:
in the first step, fix $e_1$ and prove  a result similar to Proposition \ref{step 1}.
\begin{equation}\label{step 1 hypergraph}
J(H)\lesssim_{e_i, \varepsilon, D} n_1\dots n_k \left(n_2^{-1/e_2+\varepsilon}+n_3^{-1/e_k+\varepsilon}+\dots+n_k^{-1/e_k+\varepsilon}\right)+ (n_2n_3\dots n_k)^{1+\varepsilon}.
\end{equation}
We prove this by induction by $\sum_{i=2}^k e_i$ and $\sum_{i=1}^k n_i$. For the induction step, we view the hyperedges as incidences between the grid $P_2\times\dots\times P_k\subset V_2\times\dots\times V_k\subset  \R^{d_2+\dots+d_k}$ and $n_1$ semi-algebraic sets defined by $P_1$.  If we simply apply the usual polynomial partitioning, each cell may not have the structure of a $k$-partite hypergraph. We overcome this by using  the grid polynomial partitioning: for each $i=2,\dots, k$, find a polynomial $f_i$ of degree at most $r$ to partition $P_i$ in $\R^{d_i}$ 
then take their product:
$$h(x_1,\dots, x_{d_2+\dots+d_k}):=f_2(x_1,\dots, x_{d_2}) f_3(x_{d_2+1},\dots, x_{d_2+d_3})\dots f_k(x_{d_2+\dots+{d_{k-1}}},\dots, x_{d_2+\dots+d_k})$$
By doing this, we preserve the grid structure and thus can use induction on a smaller grid in each cell.

By theorem \ref{Basu-Pollack-Roy}, for each $1<i\leq k$, $f_i$ divides $V_i$ into $O(r^{e_i})$ cells. Therefore $V_2\times\dots\times V_k \setminus Z(h)$  consists of $O(r^{e_2+\dots+e_{k}})$ cells, each cell contains a sub-grid of $P_2\times \dots\times P_{k}$ of size at most $\frac{n_2}{r^{e_2}}\times\dots\times\frac{n_{k}}{r^{e_{k}}}$. We can decompose $H$ into three parts: incidences between points in a cell and semi-algebraic sets that fully contain that cell; incidences between points in a cell and semi-algebraic sets that cross that cell; and incidences involving points lying in some $Z(f_i)$. For the first part, since those points and semi-algebraic sets form a complete subhypergraph, their representation complexity is not too large.
For the second part, since the grid structure is preserved in each cell, we can use induction assumption on smaller $\sum n_i$ to bound the representation complexity of $H$ restricted to each cell. Finally. for hyperedges involving points lying in some $Z(f_i)$, the points $P_i$ belong to $V_i\cap Z(f_i)$ which consists of several irreducible varieties of smaller dimensions, hence we can apply the induction assumption for smaller $\sum e_i$. 

From \eqref{step 1 hypergraph} we know our result holds if $n_i\geq n_1^{d_i}$ because then the term $(n_2\dots n_k)^{1+\varepsilon}$ dominates the RHS of \eqref{step 1 hypergraph}. Hence from now on we can assume $n_i\leq n_1^{d_i}$. By symmetry we can assume $n_i\leq  n_j^{d_i}$ for any distinct $i,j\leq k$. In this case the term $ E_{\vec{e}}(\vec{n})$ dominates $E_{\vec{e_I}}(\vec{n_I})\prod_{i\notin I} n_i$ for any $I\subsetneq [k]$ by Lemma \ref{E_d dominant}.
This assumption is the same with that in Remark 3.3 in \cite{Z problem hypergraph}, and that is all we need for the second step to work. 

In the second step, we fix $e_k$ and prove $J(H)\lesssim F^{*,\varepsilon}_{\vec{e}}(\vec{n})$ by induction by $e_1+\dots+e_{k-1}$ and $\sum_{i=1}^k n_i$. Again we can view the hyperedges as incidences between $n_k$ semi-algebraic sets defined by $P_k$ and the grid $P_1\times\dots\times P_{k-1}\subset \R^{d_1+\dots+d_{k-1}}$. As in the previous step, we use the grid polynomial partitioning: we can find polynomials $f_1,\dots, f_{k-1}$, each has degree at most $r$, and take their product $h=f_1\dots f_{k-1}$ so that $Z(h)$ divides $V_1\times\dots\times V_{k-1}$ into $O(r^{e_1+\dots+e_{k-1}})$ cells where each cell contains a sub-grid of $P_1\times \dots\times P_{k-1}$ of size at most $\frac{n_1}{r^{e_1}}\times\dots\times\frac{n_{k-1}}{r^{e_{k-1}}}$.

Again we can decompose $H$ into three parts. For incidences between points in a cell and semi-algebraic sets that fully contain that cell, they form complete subhypergraph and hence can be bounded easily. For incidences between points in a cell and semi-algebraic sets that crosses that cell, we apply induction for smaller $\sum n_i$ in each cell, add them and use H\"older's inequality. Here we use the fact  $E_{\vec{d}}(\vec{n})$ is the dominant term of $F^{*,\varepsilon}$, and by Lemma \ref{matrix} function $E_{\vec{d}}(\vec{n})$ behaves nicely w.r.t. the partitioning. Finally, for 
incidences involving points lying in some $Z(f_i)$,  the points $P_i$ belong to $V_i\cap Z(f_i)$ which consists of several irreducible varieties of smaller dimensions, hence we can apply the induction assumption for smaller $\sum e_i$ and use Lemma \ref{compare d and d-1}.

\end{proof}

\section{Discussion}\label{sec:conclusion}
An open question is whether the bound in theorem \ref{main thm} is tight. On one hand, when $G$ is $K_{u,u}$-free for some fixed $u$,   any lower bound on the number of edges $E(G)$ implies the same lower bound on the representation complexity $J(G)$; in particular, theorem \ref{main thm} is tight for point-line incidences (Szemer\'edi-Trotter \cite{S-T}) and not too far from tight when $Q$ is a set of hypersurfaces under certain constraints (Sheffer \cite{Sheffer-lower bound}). On the other hand, when $G$ contains too many edges, the graph may have some dense structures, such as $K_{m,n}$ or the $1/2$-degenerate point-hyperplane graphs in \cite{Elekes-Toth}, in both cases $J(G)$ is small. It might be interesting to find an example where $J(G)$ is close to its upper bound in theorem \ref{main thm} while the number of edges is much larger, or when $G$ contains some large complete bipartite subgraph.

\appendix
\section{Properties of functions $E$ and $F$}\label{property E and F}
In this appendix we present some properties of functions $E$ and $F$ defined in section \ref{sec:hypergraph}. Those properties are not difficult to prove, interested readers can find proofs in \cite{Z problem hypergraph}. Recall $E_{\vec{d}}(\vec{n})=\prod_{i=1}^k n_i^{\alpha_i}$ where  $\alpha_i=1-\frac{1/(d_i-1)}{k-1+\sum_l 1/(d_l-1)}$. 

\begin{lem}\label{matrix}
For each $i\in[k]$ we have
$\alpha_i= \sum_{j\neq i} d_j(1-\alpha_j)$. Hence the exponents $\{\alpha_i\}$ satisfy a nice system of equations:

$$
\begin{pmatrix}
1 & d_2 & \dots & d_k\\
d_1 & 1 &\dots & d_k\\
\vdots & \vdots & \ddots &\vdots\\
d_1 & d_2 & \dots & 1
\end{pmatrix}
\begin{pmatrix}
\alpha_1\\
\alpha_2\\
\vdots\\
\alpha_k
\end{pmatrix}=
\begin{pmatrix}
\sum_{i=1}^k d_i- d_1\\
\sum_{i=1}^k d_i- d_2\\
\vdots\\
\sum_{i=1}^k d_i- d_k
\end{pmatrix}
$$

As a corollary,
for any $r>0$ and each $i\in[k]$ we have 
 $$r^{d_1+\dots+d_{k-1}}E_{\vec{d}}(\frac{n_1}{r^{d_1}},\dots,\frac{n_{k-1}}{r^{d_{k-1}}},\frac{n_k}{r})=E_{\vec{d}}(\vec{n}).$$
\end{lem}

\begin{lem}\label{compare d and d-1}
Let $\mathfrak{e}_1,\dots, \mathfrak{e}_k$ be the standard basis in $\R^k$.  Then $F^\varepsilon_{\vec{d}-\mathfrak{e}_i}(\vec{n})\leq F^\varepsilon_{\vec{d}}(\vec{n})$ assuming $n_i\geq n_j^{1/(d_j)}$ for any $j\neq i$.
\end{lem}

\begin{lem}\label{E_d dominant}

 Assume $n_i< n_j^{d_i}$ for any distinct $i,j\leq k$, 
then
$ E_{\vec{d}}(\vec{n})\prod_{i=1}^k n_i^\varepsilon \geq cF^\varepsilon_{\vec{d}}(\vec{n})$ for some constant $c$. In other words, $E_{\vec{d}}(\vec{n})\prod_i n_i^\varepsilon$ is the dominant term of $F^\varepsilon_{\vec{d}}(\vec{n})$.
\end{lem}

\end{document}